\noindent\centerline{\bf Erd\'elyi-Kober Fractional Integral Operators from a Statistical Perspective -III}

\vskip.3cm\centerline{A.M. MATHAI}
\vskip.2cm\centerline{Centre for Mathematical Sciences,}
\vskip.1cm\centerline{Arunapuram P.O., Pala, Kerala-68674, India, and}
\vskip.1cm\centerline{Department of Mathematics and Statistics, McGill University,}
\vskip.1cm\centerline{Montreal, Quebec, Canada, H3A 2K6}
\vskip.2cm\centerline{and}
\vskip.2cm\centerline{H.J. HAUBOLD}
\vskip.1cm\centerline{Office for Outer Space Affairs, United Nations}
\vskip.1cm\centerline{P.O. Box 500, Vienna International Centre} 
\vskip.1cm\centerline{A - 1400 Vienna, Austria, and}
\vskip.2cm\centerline{Centre for Mathematical Sciences,}
\vskip.1cm\centerline{Arunapuram P.O., Pala, Kerala-68674, India}

\vskip.5cm\noindent{\bf Abstract}

\vskip.3cm In this article we define Kober fractional integral operators in the multivariable case. First we consider one sequence of independent random variables and an arbitrary function, which can act as the joint density of another sequence of random variables. Then we define a concept, analogous to the concept of Kober operators in the scalar variable case. This extension is achieved by using statistical techniques and the representation gives an interpretation in terms of a joint statistical density. Then we look at two sets of random variables where between the sets they are independently distributed but within each set they are dependent. Again extensions of Kober fractional integral operator are considered. Several such statistical interpretations are given for  Kober operators in the multivariable case.

\vskip.3cm\noindent{\bf 1.\hskip.3cm Introduction}

\vskip.3cm When going from a one-variable function to many-variable function there is no unique one to one correspondence. Many types of multivariable functions can be considered when one has the preselected one-variable function. Hence there is nothing called the multivariable analogue of a univariable operator. Hence we construct one multivariable operator here which is analogous to a one variable Kober fractional integral operator of the second kind. Other such analogues can be defined.

\vskip.3cm\noindent{\bf Definition 1.1.}\hskip.3cm {\it Kober fractional integral operator of the second kind in the multivariable case }\hskip.3cm This will be defined as the following fractional integral and denoted as follows:
$$\eqalignno{K_{u_j,j=1,...,k}^{(\zeta_j,\alpha_j),j=1,...,k}f(u_1,...,u_k)&=\{\prod_{j=1}^k{{u_j^{\zeta_j}}\over{\Gamma(\alpha_j)}}\}\cr
&\times\{\prod_{j=1}^k\int_{v_j=u_j}^{\infty}(v_j-u_j)^{\alpha_j-1}v_j^{-\zeta_j-\alpha_j}\}f(v_1,...,v_k){\rm d}v_1\wedge...\wedge{\rm d}v_k.&(1.1)\cr}$$

\vskip.3cm This definition is parallel to the one in the one variable case. We will now look at various connections to different problems. First we will establish a number of results in connection with statistical distribution theory. We will show that (1.1) can be treated as a constant multiple of a joint density of a number of random variables $u_1,...,u_k$ appearing in different contexts.

\vskip.3cm\noindent{\bf 1.1.\hskip.3cm Kober fractional integral operators of the second kind in multivariable case as statistical densities}

\vskip.2cm Let $x_1,x_2,...,x_k$ be independently distributed type-1 beta random variables with parameters $(\zeta_j+1,\alpha_j),j=1,...,k, \zeta_j>-1,\alpha_j>0,j=1,...,k$. Usually the parameters in a statistical density are real but the following integrals also exist for complex parameters and in that case the conditions will be $\Re(\zeta)>-1$ and $\Re(\alpha_j)>0$, $j=1,...,k$. That is, the density of $x_j$ is of the form
$$f_j(x_j)={{\Gamma(\alpha_j+\zeta_j+1)}\over{\Gamma(\alpha_j)\Gamma(\zeta_j+1)}}x_j^{\zeta_j}(1-x_j)^{\alpha_j-1},
0<x_j<1\eqno(1.2)
$$for $\alpha_j>0,\zeta_j>-1$ and $f_j(x_j)=0$ elsewhere, $j=1,...,k$ so that the joint density of $x_1,...,x_k$ is the product $f_1(x_1)...f_k(x_k)$. Let $v_1,...,v_k$ be another sequence of real scalar positive random variables having an arbitrary joint density $f(v_1,...,v_k)$. Let the two sets $(x_1,...,x_k)$ and $(v_1,...,v_k)$ be independently distributed. Consider the transformation $u_j=x_jv_j$ or $x_j={{u_j}\over{v_j}}$  and the Jacobian of the transformation is
$$\{\prod_{j=1}^k\wedge{\rm d}x_j\}\wedge\{\prod_{j=1}^k\wedge{\rm d}v_j\}=\{\prod_{j=1}^k(v_j)^{-1}\}\{\prod_{j=1}^k\wedge{\rm d}u_j\}\wedge\{\prod_{j=1}^k\wedge{\rm d}v_j\}\eqno(1.3)
$$Then the joint density of $u_1,...,u_k$, denoted by $g(u_1,...,u_k)$, is given by

$$\eqalignno{g(u_1,...,u_k)&=\{\prod_{j=1}^k{{\Gamma(\alpha_j+\zeta_j+1)}\over{\Gamma(\zeta_j+1)\Gamma(\alpha_j)}}\}
\{\prod_{j=1}^k\int_{v_j}({{u_j}\over{v_j}})^{\zeta_j}(1-{{u_j}\over{v_j}})^{\alpha_j-1}\}\cr
&\times f(v_1,...,v_k){\rm d}v_1\wedge...\wedge{\rm d}v_k\cr
&=\{\prod_{j=1}^k{{\Gamma(\alpha_j+\zeta_j+1)}\over{\Gamma(\zeta_j+1)\Gamma(\alpha_j)}}\}\{\prod_{j=1}^ku_j^{\zeta_j}
\int_{v_j=u_j}^{\infty}(v_j-u_j)^{\alpha_j-1}v_j^{-\zeta_j-\alpha_j}\}\cr
&\times f(v_1,...,v_k){\rm d}v_1\wedge...\wedge{\rm d}v_k.&(1.4)\cr}
$$Therefore one can write

\vskip.3cm\noindent{\bf Theorem 1.1.}\hskip.3cm{\it Let $x_j,u_j,v_j,j=1,...,k$ be as defined above where $x_1,...,x_k$ are independently type-1 beta distributed with parameters $(\zeta_j+1,\alpha_j),j=1,...,k$, $v_1,...,v_k$ having a joint arbitrary density $f(v_1,...,v_k)$ with $(x_1,...,x_k)$ and $(v_1,...,v_k)$ are independently distributed. If the joint density of $u_1,...,u_k$ is denoted as $g(u_1,...,u_k)$ then
$$\{\prod_{j=1}^k{{\Gamma(\zeta+1)}\over{\Gamma(\alpha_j+\zeta_j+1)}}\}g(u_1,...,u_k)
=K_{u_j,j=1,...,k}^{(\zeta_j,\alpha_j),j=1,...,k}f(u_1,...,u_k)\eqno(1.5)$$for $\Re(\zeta_j)>-1,\Re(\alpha_j)>0,j=1,...,k.$}

\vskip.3cm In this case we had $x_1,...,x_k$ mutually independently distributed and thus there were a total of $k+1$ densities involved, the $k$ of $x_1,...,x_k$ and the one of $(v_1,...,v_k)$. Let us see what happens if $x_1,...,x_k$ are not independently distributed but they have a joint density $f_1(x_1,...,x_k)$ and $(v_1,...,v_k)$ having a joint density $f_2(v_1,...,v_k)$. Then we can show that if $f_1$ can be eventually reduced to independent type-1 beta form, still we can consider Kober operators of the second kind in the multivariable case as  constant multiples of statistical densities.

\vskip.2cm Let $(x_1,...,x_k)$ have a joint type-1 Dirichlet density with parameters $(\alpha_1+1,...,\alpha_k+1;\alpha_{k+1}),\Re(\alpha_j)>-1,j=1,...,k,\Re(\alpha_{k+1})>0$, that is,

$$\eqalignno{f_1(x_1,...,x_k)&={{\Gamma(\alpha_1+...+\alpha_{k+1}+k)}
\over{\{\prod_{j=1}^k\Gamma(\alpha_j+1)\}\Gamma(\alpha_{k+1})}}x_1^{\alpha_1}...x_k^{\alpha_k}\cr
&\times (1-x_1-...-x_k)^{\alpha_{k+1}-1},0<x_j<1,j=1,...,k,0<x_1+...+x_k<1&(1.6)\cr}
$$and $f_1(x_1,...,x_k)=0$ elsewhere. Let us consider the transformations $x_1=y_1,x_2=y_2(1-y_1),...$ $x_k=(1-y_1)...(1-y_{k-1})$ or
$$\eqalignno{x_j&=y_j(1-y_1)(1-y_2)...(1-y_{j-1}),j=1,...,k\cr
\noalign{\hbox{or}}
y_j&={{x_j}\over{1-x_1-...-x_{j-1}}}, j=1,...,k.&(1.7)\cr}
$$Under this transformation the Jacobian is $(1-y_1)^{k-1}...(1-y_{k-1})$. It is easy to show that under this transformation $y_1,...,y_k$ will be independently distributed as type-1 beta variables with the parameters $(\alpha_j+1,\beta_j)$ with $\beta_j=\alpha_{j+1}+\alpha_{j+2}+...+\alpha_k+(k-j)+\alpha_{k+1}$ or $y_j$ has the density
$$f_j(y_j)={{\Gamma(\alpha_j+1+\beta_j)}\over{\Gamma(\alpha_j+1)\Gamma(\beta_j)}}y_j^{\alpha_j}(1-y_j)^{\beta_j-1},0<y_j<1
$$and $f_j(y_j)=0$ elsewhere, $\alpha_j>-1,\beta_j>0,j=1,...,k$. In the light of these observations, let us consider two sets of positive random variables $(x_1,...,x_k)$ and $(v_1,...,v_k)$ where the two sets are independently distributed with $(x_1,...,x_k)$ having a type-1 Dirichlet distribution. Let $u_j=y_jv_j=v_j({{x_j}\over{1-x_1-...-x_{j-1}}}),j=1,...,k$. Then following through the same procedure as above we have the following theorem.

\vskip.3cm\noindent{\bf Theorem 1.2.}\hskip.3cm{\it Let $(x_1,...,x_k)$ and $(v_1,...,v_k)$ be two sets of real scalar positive random variables where between sets they are independently distributed. Let $(v_1,...,v_k)$ have an arbitrary joint density $f(v_1,...,v_k)$ and let $(x_1,...,x_k)$ have a type-1 Dirichlet density with the parameters $(\alpha_1+1,...,\alpha_k+1;\alpha_{k+1})$ or with the density
$$\eqalignno{f_1(x_1,...,x_k)&=C~x_1^{\alpha_1}...x_k^{\alpha_k}(1-x_1-...-x_k)^{\alpha_{k+1}-1}&(1.8)\cr
&0<x_j<1,0<x_1+...+x_k<1,j=1,...,k\cr}
$$and $f_1(x_1,...,x_k)=0$ elsewhere, where $C$ is the normalizing constant. Let $u_j=v_j({{x_j}\over{1-x_1-...-x_{j-1}}}),j=1,...,k$. If the joint density of $u_1,...,u_k$ is again denoted by $g(u_1,...,u_k)$ then
$$\{\prod_{j=1}^k{{\Gamma(\alpha_j+1)}\over{\Gamma(\alpha_j+\beta_j+1)}}\}g(u_1,..,u_k)
=K_{u_j,j=1,...,k}^{(\alpha_j,\beta_j),j=1,...,k}f(u_1,...,u_k)\eqno(1.9)
$$where $\beta_j=\alpha_{j+1}+\alpha_{j+2}+...+\alpha_k+(k-j)+\alpha_{k+1},j=1,...,k$, $\Re(\alpha_j)>-1,j=1,...,k,\Re(\alpha_{k+1})>0,\Re(\beta_j)>0,j=1,...,k$.}

\vskip.3cm The above structure indicates that we can consider any multivariable density $f_1(x_1,...,x_k)$ for a set of  real scalar positive random variables $(x_1,...,x_k)$ and if we can find a suitable transformation to bring the joint density of the new variables as products of type-1 beta densities then the Kober fractional integral operator of the second kind for the multivariable case can be written in terms of a statistical density as shown above. There are many densities where a transformation can bring $f_1(x_1,..,x_k)$ to product of type-1 beta densities. There are several generalizations of type-1 and type-2 Dirichlet models where suitable transformations exist which can bring a set of mutually independently distributed type-1 beta random variables. We will list one more example of this type before quitting this section.

\vskip.2cm Let us consider a generalized type-1 Dirichlet model of the following type. Several types of generalizations of the following category are available.

$$\eqalignno{f_1(x_1,...,x_k)&=C_1~x_1^{\alpha_1}(1-x_1)^{\beta_1}x_2^{\alpha_2}(1-x_1-x_2)^{\beta_2}...\cr
&\times x_k^{\alpha_k}(1-x_1-...-x_k)^{\beta_k-1},0<x_1+...+x_j<1,j=1,...,k&(1.10)\cr}
$$and $f_1(x_1,...,x_k)=0$ elsewhere, where $C_1$ is a normalizing constant. Let us consider the same transformation as in (1.7). Then we can show that $y_1,...,y_k$ will be mutually independently distributed as type-1 beta random variables with the parameters $(\alpha_j+1,\gamma_j),j=1,...,k$ where $$\gamma_j=\alpha_{j+1}+\alpha_{j+2}+...+\alpha_{k-1}+\beta_j+\beta_{j+1}+...+\beta_k+(k-j)\eqno(1.11)
$$for $j=1,...,k$ with $\Re(\alpha_j)>-1,j=1,...,k$ and $\Re(\gamma_j)>0,j=1,...,k$.

\vskip.3cm\noindent{\bf Theorem 1.3.}\hskip.3cm{\it Let $x_1,...,x_k$ have a joint density of the form in (1.10). Let $v_1,..,v_k$ be another set of real scalar positive random variables having an arbitrary density $f(v_1,...,v_k)$. Between sets let $(x_1,...,x_k)$ and $(v_1,...,v_k)$ be independently distributed. Consider the transformation as in (1.7) where
$$u_j=v_j({{x_j}\over{1-x_1-...-x_{j-1}}}),j=1,...,k.
$$Let the joint density of $u_1,...,u_k$ be again denoted by $g(u_1,...,u_k)$. Then
$$\eqalignno{\{\prod{{\Gamma(\alpha_j+1)}\over{\Gamma(\alpha_j+\gamma_j+1)}}\}&g(u_1,...,u_k)\cr
&=K_{u_j,j=1,...k}^{(\alpha_j,\gamma_j),j=1,...,k}f(u_1,...,u_k)&(1.12)\cr}
$$where $\gamma_j=\alpha_{j+1}+\alpha_{j+2}+...+\alpha_k+\beta_j+\beta_{j+1}+...+\beta_k+(k-j),j=1,...,k$ for $\Re(\alpha_j)>-1,\Re(\gamma_j)>0,j=1,...,k$.}

\vskip.3cm\noindent{\bf 1.2.\hskip.3cm A Pathway Generalization of Kober Operator of the Second Kind in the Multivariable Case}

\vskip.3cm Let $x_1,...,x_k$ be independently distributed with $x_j$ having a pathway density given by
$$f_j(x_j)=c_{jp}~x_j^{\zeta_j}[1-a_j(1-q_j)x_j]^{{\eta_j}\over{1-q_j}}\eqno(1.13)
$$for $1-a_j(1-q_j)x_j>0,a_j>0,q_j<1,\eta_j>0,\zeta_j>-1$ and $f_j(x_j)=0$ elsewhere, where
$$c_{jp}={{[a_j(1-q_j)]^{\zeta_j+1}\Gamma(\zeta_j+1+{{\eta_j}\over{1-q_j}}+1)}
\over{\Gamma(\zeta_j+1)\Gamma({{\eta_j}\over{1-q_j}}+1)}}\eqno(1.14)
$$Let $v_1,...,v_k$ be real scalar positive random variables with a joint density $f(v_1,...,v_k)$. Let $(x_1,..,x_k)$ and $(v_1,...,v_k)$ be statistically independently distributed. Let $u_j=x_jv_j$, $x_j={{u_j}\over{v_j}}$, $j=1,...,k$. Then the Jacobian of the transformation is $(v_1...v_k)^{-1}$. Let the joint density of $u_1,...,u_k$ be denoted by $g(u_1,...,u_k)$. Then from the standard technique of transformation of variables the density $g$ is given by

$$\eqalignno{g(u_1,...,u_k)&=\{\prod_{j=1}^kc_{jp}u_j^{\zeta_j}\int_{v_j=a_j(1-q_j)u_j}^{\infty}v_j^{-\zeta_j-({{\eta_j}\over{1-q_j}}+1)}\cr
&\times [v_j-a(1-q_j)u_j]^{{\eta_j}\over{1-q_j}}\}f(v_1,...,v_k){\rm d}v_1\wedge...\wedge{\rm d}v_k.&(1.15)\cr}
$$Hence we may define a pathway extension of Kober operator.

\vskip.3cm\noindent{\bf Definition 1.2.}\hskip.3cm{\it A Pathway Kober Fractional Integral Operator of the Second Kind for the Multivariable Case}\hskip.3cm It will be defined and denoted as follows:
$$\eqalignno{K_{u_j,a_j,q_j,j=1,...,k}^{(\zeta_j,{{\eta_j}\over{1-q_j}}+1),j=1,...,k}f(u_1,...,u_k)&=\{\prod_{j=1}^k{{u_j^{\zeta_j}}\over{\Gamma({{\eta_j}\over{1-q_j}}+1)}}\int_{v_j>a_j(1-q_j)u_j}(v-a_j(1-q_j)u_j)^{{\eta_j}\over{1-q_j}}\}\cr
&\times f(v_1,...,v_k){\rm d}v_1\wedge...\wedge{\rm d}v_k.&(1.16)\cr}
$$for $a_j>0,q_j<1,\eta_j>0,\Re(\zeta_j)>-1$.

\vskip.3cm\noindent{\bf Theorem 1.4.}\hskip.3cm{\it Let $x_1,...,x_k$,  $v_1,...,v_k$, $u_j,j=1,...,k$ and $g(u_1,...,u_k)$ be as defined in (1.15). Let the pathway extended Kober operator be as defined in (1.16). Then

$$\eqalignno{\{\prod_{j=1}^k{{\Gamma(\zeta_j+1)}\over{\Gamma(\zeta_j+{{\eta_j}\over{1-q_j}}+2)}}\}g(u_1,...,u_k)&
=K_{u_j,a_j,q_j,j=1,...,k}^{(\zeta_j,{{\eta_j}\over{1-q_j}}+1),j=1,...,k}f(u_1,...,u_k).&(1.17)\cr}
$$}

\vskip.3cm When any particular $q_r\to 1_{-}$ then we can see the corresponding factor going to the exponential form.
$$\eqalignno{\lim_{q_r\to 1_{-}}&{{(1-q_r)^{\zeta_r+1}\Gamma(\zeta_r+{{\eta_r}\over{1-q_r}}+1)}\over{\Gamma({{\eta_r}\over{1-q_r}}+1)}}({{u_r}\over{v_r}})^{\zeta_r}{{1}\over{v_r}}
[1-a_r(1-q_r)({{u_r}\over{v_r}})]^{{\eta_r}\over{1-q_r}}\cr
&=({{u_r}\over{v_r}})^{\zeta_r}{{1}\over{v_r}}{\rm e}^{-a_r\eta_r~({{u_r}\over{v_r}})},0<v_r<\infty.&(1.18)\cr}
$$Thus, individual $q_j$'s can go to $1$ and the corresponding factor will go to exponential form or the correspondingly we get a gamma density structure for that factor.

\vskip.3cm\noindent{\bf 1.3.\hskip.3cm Mellin Transform in the Multivariable case for Kober Operators of the Second Kind}

\vskip.3cm The Mellin transform in the multivariable case is defined as
$$\eqalignno{M\{f(x_1,...,x_k);s_1,...,s_k\}&=\int_0^{\infty}...\int_0^{\infty}x_1^{s_1-1}...x_k^{s_k-1}\cr
&\times f(x_1,...,x_k){\rm d}x_1\wedge...\wedge{\rm d}x_k&(1.19)\cr}
$$whenever it exists, where $s_1,...,s_k$ in general are complex parameters. Hence for the Kober operator of the second kind we have
$$\eqalignno{M\{K_{u_j,j=1,...,k}^{(\zeta_j,\alpha_j),j=1,...,k}&f(u_1,...,u_k);s_1,...,s_k\}
=\int_0^{\infty}...\int_0^{\infty}u_1^{s_1-1}...u_k^{s_k-1}\cr
&\times\{\prod_{j=1}^k{{u_j^{\zeta_j}}\over{\Gamma(\alpha_j)}}\int_{v_j>u_j}(v_j-u_j)^{\alpha_j-1}v_j^{-\zeta-\alpha_j}\}\cr
&\times f(v_1,...,v_k){\rm d}V~{\rm d}U\cr
\noalign{\hbox{where, for example,  ${\rm d}U={\rm d}u_1\wedge...\wedge{\rm d}u_k$}}
&=\int_0^{\infty}...\int_0^{\infty}f(v_1,...,v_k)\{\prod_{j=1}^kv_j^{-\zeta_j-\alpha_j}\}[\{\prod_{j=1}^k\int_0^{v_j}u_j^{s_j+\zeta_j-1}(v_j-u_j)^{\alpha_j-1}{\rm d}u_j\}]{\rm d}V.\cr}
$$Take out $v_j$, put $y_j={{u_j}\over{v_j}}$ then the integral over $u_j$ will go to
$$v_j^{\zeta_j+\alpha_j+s_j-1}{{\Gamma(\alpha_j)\Gamma(\zeta_j+s_j)}\over{\Gamma(\alpha_j+\zeta_j+s_j)}}.$$Hence the required Mellin transform is
$$\left[\prod_{j=1}^k{{\Gamma(\zeta_j+s_j)}\over{\Gamma(\alpha_j+\zeta_j+s_j)}}\right]f^{*}(s_1,...,s_k)
$$where $f^{*}$ is the Mellin transform of $f$. Then we have the following theorem.

\vskip.3cm\noindent{\bf Theorem 1.5.}\hskip.3cm{\it For the Kober fractional integral operator of the second kind in the multivariable case as defined in (1.1) the Mellin transform, with Mellin parameters $s_1,...,s_k$ is given by
$$M\{K_{u_j,j=1,...,k}^{(\zeta_j,\alpha_j),j=1,...,k}f(u_1,...,u_k);s_1,...,s_k\}
=\left[\prod_{j=1}^k{{\Gamma(\zeta_j+s_j)}\over{\Gamma(\alpha_j+\zeta_j+s_j)}}\right]f^{*}(s_1,...,s_k)\eqno(1.20)
$$for $\Re(\alpha_j)>0,\Re(\zeta_j+s_j)>0,j=1,...,k$.}

\vskip.3cm\noindent{\bf 2.\hskip.3cm Kober Fractional Integral Operator of the First Kind for Multivariable Case}

\vskip.3cm In the multivariable case we will start with the following definition and notation.

\vskip.3cm\noindent{\bf Definition 2.1.}\hskip.3cm{\it Kober Fractional Integral Operator of the First Kind in the Multivariable Case}\hskip.3cm
$$\eqalignno{I_{u_j,j=1,...,k}^{(\zeta_j,\alpha_j),j=1,...,k}f(u_1,...,u_k)&=\{\prod_{j=1}^k{{u_j^{-\zeta_j-\alpha_j}}\over{\Gamma(\alpha_j)}}\int_{v_j=0}^{u_j}(u_j-v_j)^{\alpha_j-1}v_j^{\zeta_j}\}\cr
&\times f(v_1,...,v_k){\rm d}V.&(2.1)\cr}
$$
\vskip.3cm First, we will derive this operator as a constant multiple of a statistical density. To this end, let $x_1,...,x_k$ be independently distributed type-1 beta random variables with the parameters $(\zeta_j,\alpha_j),j=1,...,k$ or with the density
$$f_j(x_j)={{\Gamma(\zeta_j+\alpha_j)}\over{\Gamma(\zeta_j)\Gamma(\alpha_j)}}x_j^{\zeta_j-1}(1-x_j)^{\alpha_j-1},0<x_j<1,\eqno(2.2)
$$for $\alpha_j>0,\zeta_j>0$ or $\Re(\alpha_j)>0,\Re(\zeta_j)>0$ when the parameters are in the complex domain. Let $(v_1,...,v_k)$ be real scalar positive random variables having a joint density $f(v_1,...,v_k)$. Let $u_j={{v_j}\over{x_j}},j=1,...,k$. The Jacobian is given by
$${\rm d}X\wedge{\rm d}V=[\prod_{j=1}^k(-{{v_j}\over{u_j^2}})]{\rm d}U\wedge{\rm d}V\eqno(2.3)
$$where the earlier simplified notation is used. The joint density of $u_1,...,u_k$, following through the earlier steps, denoted again by $g(u_1,...,u_k)$, is given by
$$\eqalignno{g(u_1,...,u_k)&=\{\prod_{j=1}^k{{\Gamma(\zeta_j+\alpha_j)}\over{\Gamma(\zeta_j)\Gamma(\alpha_j)}}\}\{\prod_{j=1}^ku_j^{-\zeta_j-\alpha_j}\int_{v=0}^{u_j}(u_j-v_j)^{\alpha_j-1}v_j^{\zeta_j}\}\cr
&\times f(v_1,...,v_k){\rm d}V.&(2.4)\cr}
$$Hence we have the following theorem.

\vskip.3cm\noindent{\bf Theorem 2.1.}\hskip.3cm{\it Let $x_1,...,x_k$ be independently distributed type-1 beta random variables and let $v_1,...,v_k,$ $u_1,...$
$,u_k$ be as defined in (2.2) and (2.3). Let the joint density of $u_1,...,u_k$ be denoted by $g(u_1,...,u_k)$. Then
$$\{\prod_{j=1}^k{{\Gamma(\zeta_j)}\over{\Gamma(\zeta_j+\alpha_j)}}\}g(u_1,...,u_k)
=I_{u_j,j=1,...,k}^{(\zeta_j,\alpha_j),j=1,...,k}f(u_1,...,u_k).\eqno(2.5)
$$}

\vskip.3cmWe can have pathway extension to Kober operator of the first kind, parallel to the results for the case of second kind. Other properties follow parallel to those for the case of the operator of the second kind. We will evaluate the multivariable Mellin transform following through steps parallel to those in the case of the second kind and hence we give the result here as a theorem.

\vskip.3cm\noindent{\bf Theorem 2.2.}\hskip.3cm{\it The Mellin transform, with Mellin parameters $s_1,...,s_k$, for Kober fractional integral operator of the first kind in the multivariable case is given by the following:
$$M\{K_{u_j,j=1,...,k}^{(\zeta_j,\alpha_j),j=1,...,k}f(u_1,...,u_k);s_1,...,s_k\}=\{\prod_{j=1}^k{{\Gamma(1+\zeta_j-s)}\over{\Gamma(1+\alpha_j+\zeta_j-s)}}\}f^{*}(s_1,...,s_k)\eqno(2.6)
$$for $\Re(\alpha_j)>0,\Re(\zeta_j)>0,\Re(s)<1+\Re(\zeta_j),j=1,...,k$ where $f^{*}$ is the Mellin transform of $f$.}

\vskip.3cm We can obtain theorem parallel to the ones in section 1. Some of these will be stated here without proofs. The derivations are parallel to those in section 1 and hence omitted.

\vskip.3cm\noindent{\bf Theorem 2.3.}\hskip.3cm{\it Let $x_1,...,x_k$ be independently distributed as the pathway model in (1.13) with $\zeta_j$ replaced by $\zeta_j-1$, $j=1,...,k$. Let $(v_1,....,v_k)$ have a joint arbitrary density $f(v_1,...,v_k)$. Let $(x_1,...,x_k)$ and $(v_1,...,v_k)$ be independently distributed. Let $u_j={{v_j}\over{x_j}}$ or $x_j={{v_j}\over{u_j}}$,$j=1,...,k$. Let the joint density of $u_1,...,u_k$ be again denoted by $g(u_1,...,u_k)$. Then
$$\eqalignno{I_{u_j,a_j,q_j,j=1,...,k}^{(\zeta_j,{{\eta_j}\over{1-q_j}}+1),j=1,...,k}f(u_1,...,u_k)
&=\{\prod_{j=1}^k{{u_j^{-\zeta_j-({{\eta_j}\over{1-q_j}}+1)}}\over{\Gamma({{\eta_j}\over{1-q_j}}+1)}}
\int_{v_j=0}^{{{u_j}\over{1-a_j(1-q_j)}}}[u_j-a_j(1-q_j)v_j]^{{\eta_j}\over{1-q_j}}v_j^{\zeta_j}\}\cr
&\times f(v_1,...,v_k){\rm d}v_1\wedge...\wedge{\rm d}v_k&(2.8)\cr}
$$for $a_j>0,q_j<1,\eta_j>0,\Re(\zeta_j)>0,j=1,...,k$.}

\vskip.3cm From here we can have a definition for a pathway extension of Kober operator of the first kind in the multivariable case.

\vskip.3cm\noindent{\bf Definition 2.2.}\hskip.3cm{\it Kober fractional integral operator of the first kind in the multivariable case}\hskip.3cm Let the variables and parameters be as defined in Theorem 2.3. Then the pathway extended Kober operator of the first kind is defined and denoted as follows:

$$\eqalignno{I_{u_j,a_j,q_j,j=1,...,k}^{(\zeta_j,{{\eta_j}\over{1-q_j}}+1),j=1,...,k}f(u_1,...,u_k)
&=\{\prod_{j=1}^k{{u^{-\zeta_j-({{\eta_j}\over{1-q_j}}+1)}}\over{\Gamma({{\eta_j}\over{1-q_j}}+1)}}
\int_{v_j=0}^{{u_j}\over{1-a_j(1-q_j)}}[u_j-a_j(1-q_j)v_j]^{{\eta_j}\over{1-q_j}}v_j^{\zeta_j}\}\cr
&\times f(v_1,...,v_k){\rm d}v_1\wedge...\wedge{\rm d}v_k\cr}
$$for $a_j>0,q_j<1,\eta_j>0, \Re(\zeta_j)>0,j=1,...,k$.

\vskip.3cm We can list several theorems when $x_1,...,x_k$ are not independently distributed. Two such cases will be listed here without proofs. The proofs will be parallel to those in section 1 and hence omitted.

\vskip.3cm\noindent{\bf Theorem 2.4.}\hskip.3cm{\it Let $x_1,...,x_k$ have a type-1 Dirichlet density as in (1.6) with $\alpha_j$ replaced by $\alpha_j-1$ for $j=1,...,k$ and let $(v_1,...,v_k)$ have an arbitrary joint density $f(v_1,...,v_k)$ where $(x_1,...,x_k)$ and $(v_1,...,v_k)$ are independently distributed. Let $y_j={{x_j}\over{1-x_1-...-x_{j-1}}}$ and let $y_j={{v_j}\over{u_j}}, j=1,...,k$. Let the joint density of $u_1,...,u_k$ be denoted by $g(u_1,...,u_k)$. Let $\beta_j=\alpha_{j+1}+\alpha_{j+2}+...+\alpha_k$. Then
$$I_{u_j,j=1,...,k}^{(\alpha_j,\beta_j),j=1,...,k}f(u_1,...,u_k)=\{\prod_{j=1}^k{{\Gamma(\beta_j)}\over{\Gamma(\alpha_j+\beta_j)}}\}g(u_1,...,u_k).\eqno(2.9)
$$}

\vskip.3cm The next theorem is parallel to theorem 1.3 for the operator of the second kind.

\vskip.3cm\noindent{\bf Theorem 2.5.}\hskip.3cm{\it Let $x_1,...,x_k$ have a joint density as in (1.10) with $\alpha_j$ replaced by $\alpha_j-1,j=1,...,k$. Let $y_j={{x_j}\over{1-x_1-...-x_{j-1}}}$ as defined in (1.7). Let $(v_1,...,v_k)$ have an arbitrary joint density $f(v_1,...,v_k)$. Let $(x_1,...,x_k)$ and $(v_1,...,v_k)$ be independently distributed. Let $u_j={{v_j}\over{y_j}}$ or $y_j={{v_j}\over{u_j}},j=1,...,k$. Let $\gamma_j=\alpha_{j+1}+\alpha_{j+2}+...+\alpha_k+\beta_j+\beta_{j+1}+...+\beta_k$. Then
$$I_{u_j,j=1,...,k}^{(\alpha_j,\gamma_j),j=1,...,k}f(u_1,...,u_k)=\{\prod_{j=1}^k{{\Gamma(\alpha_j)}\over{\Gamma(\alpha_j+\gamma_j)}}\}g(u_1,...,u_k)\eqno(2.10)
$$for $\Re(\alpha_j)>0,\Re(\gamma_j)>0, j=1,...,k$.}

\vskip.3cm\noindent{\bf Acknowledgment}

\vskip.3cm The authors would like to thank the Department of Science and Technology, Government of India, New Delhi, for the financial assistance for this work under project number SR/S4/MS:287/05.

\vskip.3cm\centerline{\bf References}

\vskip.3cm\noindent Mathai, A.M. (1997):\hskip.3cm{\it Jacobians of Matrix Transformations and Functions of Matrix Argument}, World Scientific Publishing, New York.

\vskip.3cm\noindent Mathai, A.M. and Haubold, H.J. (2008):\hskip.3cm {\it Special Functions for Applied Scientists}, Springer, New York.

\vskip.3cm\noindent Mathai, A.M., Provost, S.B., and Hayakawa, T. (1995):\hskip.3cm{\it Bilinear Forms and Zonal Polynomials}, Springer, New York.

\bye